\documentclass{amsart}
\usepackage{todonotes}
\usepackage{mathpazo}
\usepackage{enumerate}
\usepackage{subcaption}
\usepackage{ifthen}
\usepackage{tikz,pgf,pgfmath}
	\usetikzlibrary{arrows.meta}
	\pgfdeclarelayer{background}
	\pgfdeclarelayer{middle}
	\pgfsetlayers{background,middle,main}
	\tikzstyle{enode}=[draw=black,circle,scale=\s]
	\tikzstyle{fnode}=[fill=black,draw=black,circle,scale=\s]
	\tikzstyle{tnode}=[inner sep=.2pt]
	\tikzstyle{edge}=[line width=1pt]
	\tikzstyle{label}=[fill=white,text=\red,circle,inner sep=.2pt,scale=.67]

	\newcommand{\red}{red!50!gray}
	
\usepackage[colorlinks=true,
	urlcolor=blue!20!black,
	citecolor=green!50!black,
	pagebackref=true]{hyperref}
\usepackage{amssymb}
\usepackage[nobysame]{amsrefs}
\newtheorem{theorem}{Theorem}[section]
\newtheorem{proposition}[theorem]{Proposition}

\newtheorem{lemma}[theorem]{Lemma}
\newtheorem{corollary}[theorem]{Corollary}
\newtheorem{question}[theorem]{Question}
\theoremstyle{remark}

\newcommand{\defn}[1]{{\color{green!50!black}\emph{#1}}}
\newcommand{\defs}{\stackrel{\mathsf{def}}{=}}
\newcommand{\ie}{i.e.\;}
\newcommand{\Lattice}{\textbf{L}}
\newcommand{\dual}{\mathsf{d}}
\newcommand{\len}{\mathsf{len}}
\newcommand{\least}{\hat{0}}
\newcommand{\grtst}{\hat{1}}
\newcommand{\Covers}{\mathsf{Cov}}
\newcommand{\JI}{\mathsf{J}}
\newcommand{\MI}{\mathsf{M}}
\newcommand{\perspective}{\doublebarwedge}
\newcommand{\gammairr}{\gamma}
\newcommand{\lmlabeling}[1]{\lambda_{#1}}
\newcommand{\horArrow}[8]{
	\coordinate(start) at (#1,#4);
	\coordinate(middle) at (#2,#4);
	\coordinate(end) at (#3,#4);
	\draw[#5,line width=.75pt,double distance=1.5pt,arrows={-Latex[fill=white,length=1pt 2 0]}] (start) -- (end);
	\ifthenelse{\equal{#6}{1}}{
		\draw[#5](middle) -- +(-.2*\x,-.1*\y) -- +(.2*\x,.1*\y);
	}{}
	\ifthenelse{\equal{#8}{1}}{
		\draw(middle)++(0,.15*\y) node[scale=.5]{#7};
	}{
		\draw(middle)++(0,-.15*\y) node[scale=.5]{#7};
	}
}

\newcommand{\verArrow}[8]{
	\coordinate(start) at (#1,#2);
	\coordinate(middle) at (#1,#3);
	\coordinate(end) at (#1,#4);
	\draw[#5,line width=.75pt,double distance=1.5pt,arrows={-Latex[fill=white,length=1pt 2 0]}] (start) -- (end);
	\ifthenelse{\equal{#6}{1}}{
		\draw[#5](middle) -- +(-.1*\x,-.1*\y) -- +(.1*\x,.1*\y);
	}{}
	\ifthenelse{\equal{#8}{1}}{
		\draw(middle)++(-.15*\x,0) node[rotate=90,scale=.5]{#7};
	}{
		\draw(middle)++(.15*\x,0) node[rotate=90,scale=.5]{#7};
	}
}

\newcommand{\poset}[5]{
	\begin{tikzpicture}
		\def\d{#1};
		\ifthenelse{\equal{#3}{1}}{
			\filldraw[draw=black,fill=#4](.65*\d,.75*\d) -- (3.35*\d,.75*\d) -- (3.35*\d,2.25*\d) -- (.65*\d,2.25*\d) -- cycle;
		}{
			\fill[white](.65*\d,.75*\d) -- (3.35*\d,.75*\d) -- (3.35*\d,2.25*\d) -- (.65*\d,2.25*\d) -- cycle;
		}
		\coordinate(p1) at (1*\d,1*\d);
		\coordinate(p2) at (2*\d,1*\d);
		\coordinate(p3) at (3*\d,1*\d);
		\coordinate(p4) at (1.5*\d,2*\d);
		\draw[thick](p1) -- (p4);
		\draw[thick](p2) -- (p4);
		\ifthenelse{\equal{#5}{1}}{
			\draw(p1) node[draw=black,fill=white,circle,scale=1*\d]{$1$};
			\draw(p2) node[draw=black,fill=white,circle,scale=1*\d]{$2$};
			\draw(p3) node[draw=black,fill=white,circle,scale=1*\d]{$3$};
			\draw(p4) node[draw=black,fill=white,circle,scale=1*\d]{$4$};
			
		}{
			\draw(p1) node[draw=black,fill=white,circle,scale=1*\d]{};
			\draw(p2) node[draw=black,fill=white,circle,scale=1*\d]{};
			\draw(p3) node[draw=black,fill=white,circle,scale=1*\d]{};
			\draw(p4) node[draw=black,fill=white,circle,scale=1*\d]{};
		}
		\foreach \a in {#2}{
			\draw(p\a) node[draw=black,fill=black,circle,scale=1*\d]{};
		}
	\end{tikzpicture}
}

\title{Extremality, Left-Modularity and Semidistributivity}
\author{Henri M{\"u}hle}
\address{Technische Universit{\"a}t Dresden, Institut f{\"u}r Algebra, Zellescher Weg 12--14, 01069 Dresden, Germany.}
\email{henri.muehle@tu-dresden.de}

\begin{document}

\begin{abstract}
	In this article we investigate the relations between three classes of lattices each extending the class of distributive lattices in a different way.  In particular, we consider join-semidistributive, join-extremal and left-modular lattices, respectively.  Our main motivation is a recent result by Thomas and Williams proving that every semidistributive, extremal lattice is left modular.  We prove the converse of this on a slightly more general level.  Our main result asserts that every join-semidistributive, left-modular lattice is join extremal.  We also relate these properties to the topological notion of lexicographic shellability.  
\end{abstract}

\maketitle

\section{Introduction}
    \label{sec:introduction}
One of the most fundamental (and at the same time most important) classes of lattices is the class of distributive lattices.  These lattices are characterized by satisfying the well-known distributive laws with respect to meet and join operations.  More precisely, a lattice is \emph{distributive} if for any choice of elements $a,b,c$ it holds that
\begin{align}
    (a\vee b)\wedge(a\vee c) = a\vee(b\wedge c),\\
    (a\wedge b)\vee(a\wedge c) = a\wedge(b\vee c).
\end{align}
A celebrated result of Birkhoff's states that a finite lattice is distributive if and only if it is isomorphic to the lattice of order ideals of some (\emph{underlying}) finite partially ordered set (or \emph{poset})~\cite{birkhoff37rings}.

\bigskip

From this perspective it is straightforward to describe the covering pairs in a distributive lattice.  In particular, an order ideal $a$ is covered by some other order ideal $b$ if and only if $a$ is obtained by removing a maximal element from $b$.  We can therefore realize $a$ as the join of the order ideals generated by the maximal elements of $a$.
This representation is \emph{canonical} in the sense that it is minimal in size and contains elements as close to the bottom as possible.  Of course, this idea of representing lattice elements in terms of as few as possible elements which are as far down as possible, may be considered for arbitrary lattices.  Such \emph{canonical join representations} play an important role in the solution of the word problem for free lattices~\cites{whitman41free,whitman42free}.  There is even a nice characterization of the finite lattices in which every element admits such a canonical join representation.  These are precisely the finite \emph{join semidistributive} lattices, \ie lattices where every three elements $a,b,c$ satisfy the following implication:
\begin{align}
    \text{if}\;a\vee b=a\vee c\quad\text{then}\quad a\vee b=a\vee(b\wedge c).
\end{align}
See \cite{freese95free}*{Chapter~II.5} for more background.

\bigskip

The order-ideal-representation of distributive lattices exhibits yet another intriguing property.  By the reasoning from the previous paragraph, we may conclude that the order ideals generated by single elements of the underlying poset are \emph{join irreducible}, \ie they cannot be expressed as a nontrivial join of other elements.  In fact, all join-irreducible order ideals are of this form.  Thus, the number of join-irreducible elements in a distributive lattice agrees with the size of the underlying poset.  In a related fashion, we can use any linear extension of the underlying poset to construct a maximal chain in a distributive lattice.  Since for distributive lattices, all maximal chains have the same size, this implies that the maximum size of a maximal chain equals the number of join-irreducible elements.  Following \cite{markowsky92primes}, finite lattices with this property are called \emph{join extremal}.

\bigskip

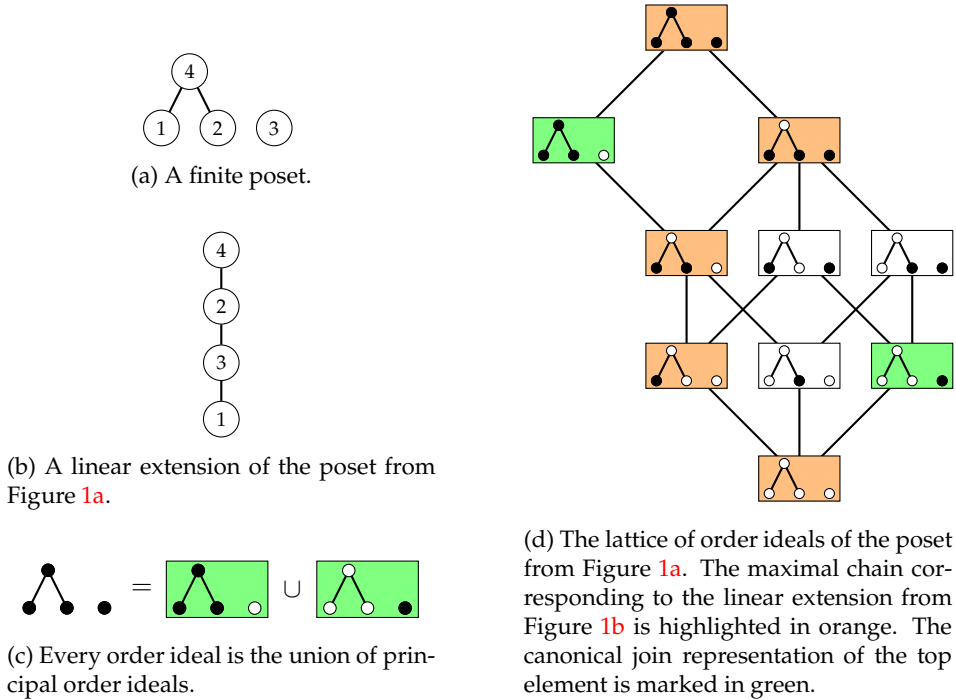
\begin{figure}
	\centering
	\begin{minipage}[b]{1\textwidth}
		\begin{minipage}[b]{.45\textwidth}
			\begin{subfigure}[t]{1\textwidth}
				\centering
				\poset{.75}{}{0}{white}{1}
				\caption{A finite poset.}
				\label{fig:poset}
			\end{subfigure}

			\vspace*{.5cm}

			\begin{subfigure}[b]{1\textwidth}
				\centering
				\begin{tikzpicture}
					\draw(.75,.75) node[draw,circle,scale=.75](n1){$1$};
					\draw(.75,1.5) node[draw,circle,scale=.75](n2){$3$};
					\draw(.75,2.25) node[draw,circle,scale=.75](n3){$2$};
					\draw(.75,3) node[draw,circle,scale=.75](n4){$4$};
					\draw[thick](n1) -- (n2) -- (n3) -- (n4);
				\end{tikzpicture}
				\caption{A linear extension of the poset from Figure~\ref{fig:poset}.}
				\label{fig:linear_extension}
			\end{subfigure}

			\vspace*{.5cm}

			\begin{subfigure}[b]{1\textwidth}
				\begin{tikzpicture}
					\def\x{2};
					\draw(1*\x,1*\x) node{\poset{.5}{1,2,3,4}{0}{white}{0}};
					\draw(2*\x,1*\x) node{\poset{.5}{1,2,4}{1}{white!50!green}{0}};
					\draw(3*\x,1*\x) node{\poset{.5}{3}{1}{white!50!green}{0}};
					\draw(1.5*\x,1*\x) node{$=$};
					\draw(2.5*\x,1*\x) node{$\cup$};
				\end{tikzpicture}
				\caption{Every order ideal is the union of principal order ideals.}
				\label{fig:union}
			\end{subfigure}
		\end{minipage}
		\hspace*{1cm}
		\begin{minipage}[b]{.45\textwidth}
			\begin{subfigure}[b]{1\textwidth}
				\centering
				\begin{tikzpicture}
					\def\x{1.5};
					\def\s{.4};
					\coordinate(n1) at (2*\x,1*\x);
					\coordinate(n2) at (1*\x,2*\x);
					\coordinate(n3) at (2*\x,2*\x);
					\coordinate(n4) at (3*\x,2*\x);
					\coordinate(n5) at (1*\x,3*\x);
					\coordinate(n6) at (2*\x,3*\x);
					\coordinate(n7) at (3*\x,3*\x);
					\coordinate(n8) at (0*\x,4*\x);
					\coordinate(n9) at (2*\x,4*\x);
					\coordinate(n10) at (1*\x,5*\x);
					\draw[thick](n1) -- (n2);
					\draw[thick](n1) -- (n3);
					\draw[thick](n1) -- (n4);
					\draw[thick](n2) -- (n5);
					\draw[thick](n2) -- (n6);
					\draw[thick](n3) -- (n5);
					\draw[thick](n3) -- (n7);
					\draw[thick](n4) -- (n6);
					\draw[thick](n4) -- (n7);
					\draw[thick](n5) -- (n8);
					\draw[thick](n5) -- (n9);
					\draw[thick](n6) -- (n9);
					\draw[thick](n7) -- (n9);
					\draw[thick](n8) -- (n10);
					\draw[thick](n9) -- (n10);
					\draw(n1) node{\poset{\s}{}{1}{white!50!orange}{0}};
					\draw(n2) node{\poset{\s}{1}{1}{white!50!orange}{0}};
					\draw(n3) node{\poset{\s}{2}{1}{white}{0}};
					\draw(n4) node{\poset{\s}{3}{1}{white!50!green}{0}};
					\draw(n5) node{\poset{\s}{1,2}{1}{white!50!orange}{0}};
					\draw(n6) node{\poset{\s}{1,3}{1}{white}{0}};
					\draw(n7) node{\poset{\s}{2,3}{1}{white}{0}};
					\draw(n8) node{\poset{\s}{1,2,4}{1}{white!50!green}{0}};
					\draw(n9) node{\poset{\s}{1,2,3}{1}{white!50!orange}{0}};
					\draw(n10) node{\poset{\s}{1,2,3,4}{1}{white!50!orange}{0}};
				\end{tikzpicture}
				\caption{The lattice of order ideals of the poset from Figure~\ref{fig:poset}.  The maximal chain corresponding to the linear extension from Figure~\ref{fig:linear_extension} is highlighted in orange.  The canonical join representation of the top element is marked in green.}
				\label{fig:lattice}
			\end{subfigure}
		\end{minipage}
	\end{minipage}
	\caption{Illustrating several aspects of distributive lattices.}
	\label{fig:distributive_overview}
\end{figure}

Another remarkable property of distributive lattices is the fact that any three elements $a,b,c$, where $b<c$, satisfy the \emph{modular} equality:
\begin{align}\label{eq:modular}
    (b\vee a)\wedge c = b\vee(a\wedge c).
\end{align}
The element $a$ is then called \emph{left modular}.  Drawing inspiration from group theory, Stanley introduced the class of \emph{supersolvable} lattices~\cite{stanley72supersolvable}.  These are graded lattices which possess a maximal chain consisting entirely of left-modular elements.  In particular, every distributive lattice is supersolvable.  Several researchers have subsequently studied non-graded lattices which have a maximal chain of left-modular elements, see for instance \cites{liu00left,mcnamara06poset,thomas05graded}.  Such lattices are themselves called \emph{left modular} and have several intriguing combinatorial and topological properties.

\bigskip

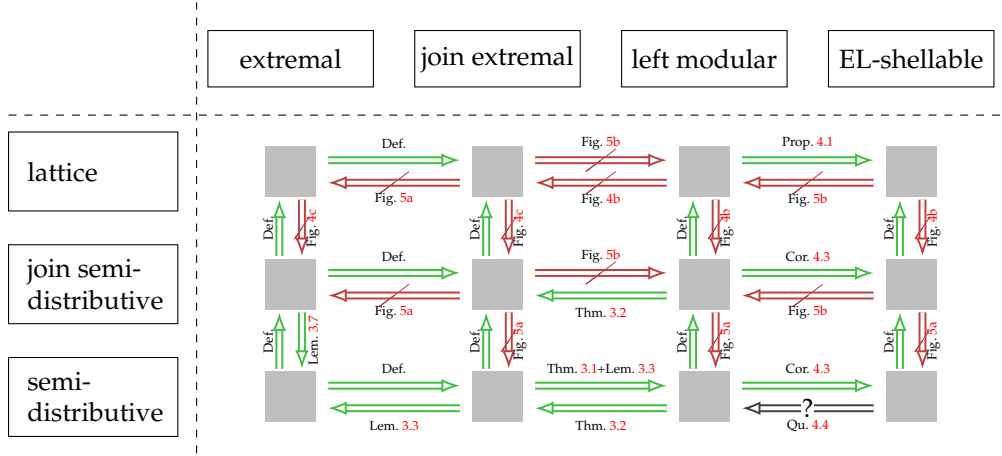
\begin{figure}
	\centering
	\begin{tikzpicture}
		\def\x{1};
		\def\y{1.5};
		\def\dx{2.75};
		\def\dy{1};
		\draw({(.75+\dx)*\x},{(1+4*\dy)*\y}) node{extremal};
			\draw({(-.35+\dx)*\x},{(.75+4*\dy)*\y}) -- ({(1.85+\dx)*\x},{(.75+4*\dy)*\y}) -- ({(1.85+\dx)*\x},{(1.25+4*\dy)*\y}) -- ({(-.35+\dx)*\x},{(1.25+4*\dy)*\y}) -- cycle;
		\draw({(.75+2*\dx)*\x},{(1+4*\dy)*\y}) node{join extremal};
			\draw({(-.35+2*\dx)*\x},{(.75+4*\dy)*\y}) -- ({(1.85+2*\dx)*\x},{(.75+4*\dy)*\y}) -- ({(1.85+2*\dx)*\x},{(1.25+4*\dy)*\y}) -- ({(-.35+2*\dx)*\x},{(1.25+4*\dy)*\y}) -- cycle;
		\draw({(.75+3*\dx)*\x},{(1+4*\dy)*\y}) node{left modular};
			\draw({(-.35+3*\dx)*\x},{(.75+4*\dy)*\y}) -- ({(1.85+3*\dx)*\x},{(.75+4*\dy)*\y}) -- ({(1.85+3*\dx)*\x},{(1.25+4*\dy)*\y}) -- ({(-.35+3*\dx)*\x},{(1.25+4*\dy)*\y}) -- cycle;
		\draw({(.75+4*\dx)*\x},{(1+4*\dy)*\y}) node{EL-shellable};
			\draw({(-.35+4*\dx)*\x},{(.75+4*\dy)*\y}) -- ({(1.85+4*\dx)*\x},{(.75+4*\dy)*\y}) -- ({(1.85+4*\dx)*\x},{(1.25+4*\dy)*\y}) -- ({(-.35+4*\dx)*\x},{(1.25+4*\dy)*\y}) -- cycle;
		\draw(1*\x,{(1+3*\dy)*\y}) node{\parbox{2cm}{lattice}};
			\draw(-.25*\x,{(.65+3*\dy)*\y}) -- (2*\x,{(.65+3*\dy)*\y}) -- (2*\x,{(1.35+3*\dy)*\y}) -- (-.25*\x,{(1.35+3*\dy)*\y}) -- cycle;
		\draw(1*\x,{(1+2*\dy)*\y}) node{\parbox{2cm}{join semi-\\distributive}};
			\draw(-.25*\x,{(.65+2*\dy)*\y}) -- (2*\x,{(.65+2*\dy)*\y}) -- (2*\x,{(1.35+2*\dy)*\y}) -- (-.25*\x,{(1.35+2*\dy)*\y}) -- cycle;
		\draw(1*\x,{(1+1*\dy)*\y}) node{\parbox{2cm}{semi-\\distributive}};
			\draw(-.25*\x,{(.65+1*\dy)*\y}) -- (2*\x,{(.65+1*\dy)*\y}) -- (2*\x,{(1.35+1*\dy)*\y}) -- (-.25*\x,{(1.35+1*\dy)*\y}) -- cycle;
		\draw[dashed](-.25*\x,{(.5+4*\dy)*\y}) -- ({(2+4*\dx)*\x},{(.5+4*\dy)*\y});
		\draw[dashed](2.25*\x,1.5*\y) -- (2.25*\x,{(.5+5*\dy)*\y});
		\coordinate(a11) at ({.75+\dx)*\x},{(1+3*\dy)*\y});
		\coordinate(a12) at ({.75+2*\dx)*\x},{(1+3*\dy)*\y});
		\coordinate(a13) at ({.75+3*\dx)*\x},{(1+3*\dy)*\y});
		\coordinate(a14) at ({.75+4*\dx)*\x},{(1+3*\dy)*\y});
		\coordinate(a21) at ({.75+\dx)*\x},{(1+2*\dy)*\y});
		\coordinate(a22) at ({.75+2*\dx)*\x},{(1+2*\dy)*\y});
		\coordinate(a23) at ({.75+3*\dx)*\x},{(1+2*\dy)*\y});
		\coordinate(a24) at ({.75+4*\dx)*\x},{(1+2*\dy)*\y});
		\coordinate(a31) at ({.75+\dx)*\x},{(1+1*\dy)*\y});
		\coordinate(a32) at ({.75+2*\dx)*\x},{(1+1*\dy)*\y});
		\coordinate(a33) at ({.75+3*\dx)*\x},{(1+1*\dy)*\y});
		\coordinate(a34) at ({.75+4*\dx)*\x},{(1+1*\dy)*\y});
		\draw(a11) node[fill=gray!50!white,rectangle,minimum width=.67cm,minimum height=.67cm]{};
		\draw(a12) node[fill=gray!50!white,rectangle,minimum width=.67cm,minimum height=.67cm]{};
		\draw(a13) node[fill=gray!50!white,rectangle,minimum width=.67cm,minimum height=.67cm]{};
		\draw(a14) node[fill=gray!50!white,rectangle,minimum width=.67cm,minimum height=.67cm]{};
		\draw(a21) node[fill=gray!50!white,rectangle,minimum width=.67cm,minimum height=.67cm]{};
		\draw(a22) node[fill=gray!50!white,rectangle,minimum width=.67cm,minimum height=.67cm]{};
		\draw(a23) node[fill=gray!50!white,rectangle,minimum width=.67cm,minimum height=.67cm]{};
		\draw(a24) node[fill=gray!50!white,rectangle,minimum width=.67cm,minimum height=.67cm]{};
		\draw(a31) node[fill=gray!50!white,rectangle,minimum width=.67cm,minimum height=.67cm]{};
		\draw(a32) node[fill=gray!50!white,rectangle,minimum width=.67cm,minimum height=.67cm]{};
		\draw(a33) node[fill=gray!50!white,rectangle,minimum width=.67cm,minimum height=.67cm]{};
		\draw(a34) node[fill=gray!50!white,rectangle,minimum width=.67cm,minimum height=.67cm]{};
		\horArrow{{1.25+1*\dx)*\x}}{{.75+1.5*\dx)*\x}}{{.25+2*\dx)*\x}}{{(1.1+3*\dy)*\y}}{green!50!gray}{0}{Def.}{1}
		\horArrow{{.25+2*\dx)*\x}}{{.75+1.5*\dx)*\x}}{{1.25+1*\dx)*\x}}{{(.9+3*\dy)*\y}}{red!50!gray}{1}{Fig.~\ref{fig:left_modular_not_semidistributive}}{0}
		\horArrow{{1.25+2*\dx)*\x}}{{.75+2.5*\dx)*\x}}{{.25+3*\dx)*\x}}{{(1.1+3*\dy)*\y}}{red!50!gray}{1}{Fig.~\ref{fig:join_extremal_not_left_modular}}{1}
		\horArrow{{.25+3*\dx)*\x}}{{.75+2.5*\dx)*\x}}{{1.25+2*\dx)*\x}}{{(.9+3*\dy)*\y}}{red!50!gray}{1}{Fig.~\ref{fig:m3}}{0}
		\horArrow{{1.25+3*\dx)*\x}}{{.75+3.5*\dx)*\x}}{{.25+4*\dx)*\x}}{{(1.1+3*\dy)*\y}}{green!50!gray}{0}{Prop.~\ref{prop:left_modular_shellable}}{1}
		\horArrow{{.25+4*\dx)*\x}}{{.75+3.5*\dx)*\x}}{{1.25+3*\dx)*\x}}{{(.9+3*\dy)*\y}}{red!50!gray}{1}{Fig.~\ref{fig:join_extremal_not_left_modular}}{0}
		\horArrow{{1.25+1*\dx)*\x}}{{.75+1.5*\dx)*\x}}{{.25+2*\dx)*\x}}{{(1.1+2*\dy)*\y}}{green!50!gray}{0}{Def.}{1}
		\horArrow{{.25+2*\dx)*\x}}{{.75+1.5*\dx)*\x}}{{1.25+1*\dx)*\x}}{{(.9+2*\dy)*\y}}{red!50!gray}{1}{Fig.~\ref{fig:left_modular_not_semidistributive}}{0}
		\horArrow{{1.25+2*\dx)*\x}}{{.75+2.5*\dx)*\x}}{{.25+3*\dx)*\x}}{{(1.1+2*\dy)*\y}}{red!50!gray}{1}{Fig.~\ref{fig:join_extremal_not_left_modular}}{1}
		\horArrow{{.25+3*\dx)*\x}}{{.75+2.5*\dx)*\x}}{{1.25+2*\dx)*\x}}{{(.9+2*\dy)*\y}}{green!50!gray}{0}{Thm.~\ref{thm:semidistributive_left_modular_is_trim}}{0}
		\horArrow{{1.25+3*\dx)*\x}}{{.75+3.5*\dx)*\x}}{{.25+4*\dx)*\x}}{{(1.1+2*\dy)*\y}}{green!50!gray}{0}{Cor.~\ref{cor:left_modular_shellable}}{1}
		\horArrow{{.25+4*\dx)*\x}}{{.75+3.5*\dx)*\x}}{{1.25+3*\dx)*\x}}{{(.9+2*\dy)*\y}}{red!50!gray}{1}{Fig.~\ref{fig:join_extremal_not_left_modular}}{0}
		\horArrow{{1.25+1*\dx)*\x}}{{.75+1.5*\dx)*\x}}{{.25+2*\dx)*\x}}{{(1.1+1*\dy)*\y}}{green!50!gray}{0}{Def.}{1}
		\horArrow{{.25+2*\dx)*\x}}{{.75+1.5*\dx)*\x}}{{1.25+1*\dx)*\x}}{{(.9+1*\dy)*\y}}{green!50!gray}{0}{Lem.~\ref{lem:semidistributive_extremal}}{0}
		\horArrow{{1.25+2*\dx)*\x}}{{.75+2.5*\dx)*\x}}{{.25+3*\dx)*\x}}{{(1.1+1*\dy)*\y}}{green!50!gray}{0}{Thm.~\ref{thm:extremal_is_trim}+Lem.~\ref{lem:semidistributive_extremal}}{1}
		\horArrow{{.25+3*\dx)*\x}}{{.75+2.5*\dx)*\x}}{{1.25+2*\dx)*\x}}{{(.9+1*\dy)*\y}}{green!50!gray}{0}{Thm.~\ref{thm:semidistributive_left_modular_is_trim}}{0}
		\horArrow{{1.25+3*\dx)*\x}}{{.75+3.5*\dx)*\x}}{{.25+4*\dx)*\x}}{{(1.1+1*\dy)*\y}}{green!50!gray}{0}{Cor.~\ref{cor:left_modular_shellable}}{1}
		\horArrow{{.25+4*\dx)*\x}}{{.75+3.5*\dx)*\x}}{{1.25+3*\dx)*\x}}{{(.9+1*\dy)*\y}}{black!50!gray}{0}{Qu.~\ref{qu:shellable_left_modular}}{0}
			\draw({.75+3.5*\dx)*\x},{(.9+1*\dy)*\y}) node[fill=white,inner sep=.2pt]{?};
		\verArrow{{.6+1*\dx)*\x}}{{(1.25+2*\dy)*\y}}{{(1+2.5*\dy)*\y}}{{(.75+3*\dy)*\y}}{green!50!gray}{0}{Def.}{1}
		\verArrow{{.9+1*\dx)*\x}}{{(.75+3*\dy)*\y}}{{(1+2.5*\dy)*\y}}{{(1.25+2*\dy)*\y}}{red!50!gray}{1}{Fig.~\ref{fig:extremal}}{0}
		\verArrow{{.6+2*\dx)*\x}}{{(1.25+2*\dy)*\y}}{{(1+2.5*\dy)*\y}}{{(.75+3*\dy)*\y}}{green!50!gray}{0}{Def.}{1}
		\verArrow{{.9+2*\dx)*\x}}{{(.75+3*\dy)*\y}}{{(1+2.5*\dy)*\y}}{{(1.25+2*\dy)*\y}}{red!50!gray}{1}{Fig.~\ref{fig:extremal}}{0}
		\verArrow{{.6+3*\dx)*\x}}{{(1.25+2*\dy)*\y}}{{(1+2.5*\dy)*\y}}{{(.75+3*\dy)*\y}}{green!50!gray}{0}{Def.}{1}
		\verArrow{{.9+3*\dx)*\x}}{{(.75+3*\dy)*\y}}{{(1+2.5*\dy)*\y}}{{(1.25+2*\dy)*\y}}{red!50!gray}{1}{Fig.~\ref{fig:m3}}{0}
		\verArrow{{.6+4*\dx)*\x}}{{(1.25+2*\dy)*\y}}{{(1+2.5*\dy)*\y}}{{(.75+3*\dy)*\y}}{green!50!gray}{0}{Def.}{1}
		\verArrow{{.9+4*\dx)*\x}}{{(.75+3*\dy)*\y}}{{(1+2.5*\dy)*\y}}{{(1.25+2*\dy)*\y}}{red!50!gray}{1}{Fig.~\ref{fig:m3}}{0}
		\verArrow{{.6+1*\dx)*\x}}{{(1.25+1*\dy)*\y}}{{(1+1.5*\dy)*\y}}{{(.75+2*\dy)*\y}}{green!50!gray}{0}{Def.}{1}
		\verArrow{{.9+1*\dx)*\x}}{{(.75+2*\dy)*\y}}{{(1+1.5*\dy)*\y}}{{(1.25+1*\dy)*\y}}{green!50!gray}{0}{Lem.~\ref{lem:extremal_semidistributive}}{0}
		\verArrow{{.6+2*\dx)*\x}}{{(1.25+1*\dy)*\y}}{{(1+1.5*\dy)*\y}}{{(.75+2*\dy)*\y}}{green!50!gray}{0}{Def.}{1}
		\verArrow{{.9+2*\dx)*\x}}{{(.75+2*\dy)*\y}}{{(1+1.5*\dy)*\y}}{{(1.25+1*\dy)*\y}}{red!50!gray}{1}{Fig.~\ref{fig:left_modular_not_semidistributive}}{0}
		\verArrow{{.6+3*\dx)*\x}}{{(1.25+1*\dy)*\y}}{{(1+1.5*\dy)*\y}}{{(.75+2*\dy)*\y}}{green!50!gray}{0}{Def.}{1}
		\verArrow{{.9+3*\dx)*\x}}{{(.75+2*\dy)*\y}}{{(1+1.5*\dy)*\y}}{{(1.25+1*\dy)*\y}}{red!50!gray}{1}{Fig.~\ref{fig:left_modular_not_semidistributive}}{0}
		\verArrow{{.6+4*\dx)*\x}}{{(1.25+1*\dy)*\y}}{{(1+1.5*\dy)*\y}}{{(.75+2*\dy)*\y}}{green!50!gray}{0}{Def.}{1}
		\verArrow{{.9+4*\dx)*\x}}{{(.75+2*\dy)*\y}}{{(1+1.5*\dy)*\y}}{{(1.25+1*\dy)*\y}}{red!50!gray}{1}{Fig.~\ref{fig:left_modular_not_semidistributive}}{0}
	\end{tikzpicture}
	\caption{Implications among the various types of lattices.}
	\label{fig:big_picture}
\end{figure}

The three properties of distributive lattices that we have just reviewed are illustrated in Figure~\ref{fig:distributive_overview}.  The purpose of this article is the study of the interactions between the three induced lattice classes: join-semidistributive lattices, join-extremal lattices and left-modular lattices.  As described above, each of these families contains the class of distributive lattices, but none of them is contained in another.  The main motivation for the research presented here is the recent result which states that every semidistributive extremal lattice is left modular~\cite{thomas19rowmotion}*{Theorem~1.4}.  Our main result is the converse of this claim on the level of join-semidistributive lattices: in Theorem~\ref{thm:semidistributive_left_modular_is_trim} we prove that every join-semidistributive, left modular lattice is join-extremal.  
An overview of the relations between various combinations of the considered properties is given in Figure~\ref{fig:big_picture}.  Finally, we consider the interaction of these properties with the topological property of lexicographic shellability.  We end with an open question asking whether any semidistributive, lexicographically shellable lattice is necessarily left modular (Question~\ref{qu:shellable_left_modular}).

\section{Lattice-theoretic preliminaries}
	\label{sec:lattice_theory}
In this note we consider only finite lattices, and we refer the interested reader to \cites{davey02introduction,freese95free,gratzer78general} for further background information.  Moreover, we usually view lattices from an order-theoretic perspective.  More precisely, if $\Lattice=(L;\vee,\wedge,\least,\grtst)$ is a (finite) lattice given as an algebraic structure with signature $(2,2,0,0)$, then we may consider this as a partially ordered set by ordering its elements via $a\leq b$ if and only if $a\vee b=b$ (or equivalently $a\wedge b=a$).  In this case, we write $\Lattice=(L,\leq)$ instead.
%
The \defn{dual} lattice is $\Lattice^{\dual}\defs(L,\geq)$.  For an integer $n>0$, we define $[n]\defs\{1,2,\ldots,n\}$.

\subsection{Covering pairs, perspectivity and join-irreducibles}
	\label{sec:perspectivity}
Two elements $a,b\in L$ form a \defn{covering pair} if $a<b$ and there does not exist $c\in L$ such that $a<c<b$.  In that event we write $a\lessdot b$.  Then, $a$ \defn{is covered by} $b$ and $b$ \defn{covers} $a$.  The \defn{cover relation} of $\Lattice$ is the set of covering pairs, defined by
\begin{displaymath}
	\Covers(\Lattice) \defs \bigl\{(a,b)\colon a\lessdot b\bigr\} \subseteq L\times L.
\end{displaymath}

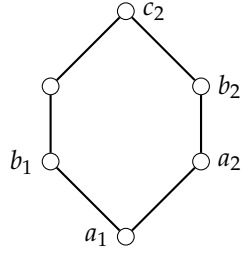
\begin{figure}
	\centering
	\begin{tikzpicture}
		\def\x{1};
		\draw(2*\x,1*\x) node[draw,circle,scale=.67](a1){};
			\draw(1.9*\x,1*\x) node[anchor=east]{$a_{1}$};
		\draw(1*\x,2*\x) node[draw,circle,scale=.67](a2){};
			\draw(.9*\x,2*\x) node[anchor=east]{$b_{1}$};
		\draw(3*\x,2*\x) node[draw,circle,scale=.67](a3){};
			\draw(3.1*\x,2*\x) node[anchor=west]{$a_{2}$};
		\draw(1*\x,3*\x) node[draw,circle,scale=.67](a4){};
		\draw(3*\x,3*\x) node[draw,circle,scale=.67](a5){};
			\draw(3.1*\x,3*\x) node[anchor=west]{$b_{2}$};
		\draw(2*\x,4*\x) node[draw,circle,scale=.67](a6){};
			\draw(2.1*\x,4*\x) node[anchor=west]{$c_{2}$};
		\draw[thick](a1) -- (a2);
		\draw[thick](a1) -- (a3);
		\draw[thick](a2) -- (a4);
		\draw[thick](a3) -- (a5);
		\draw[thick](a4) -- (a6);
		\draw[thick](a5) -- (a6);
	\end{tikzpicture}
	\caption{The covering pairs $(a_{1},b_{1})$ and $(b_{2},c_{2})$ are perspective.  The covering pair $(a_{2},b_{2})$ is not perspective to any other covering pair.}
	\label{fig:perspective_covers}
\end{figure}

Two covering pairs $(a_{1},b_{1})$ and $(a_{2},b_{2})$ are \defn{perspective} if either $b_{1}\vee a_{2}=b_{2}$ and $b_{1}\wedge a_{2}=a_{1}$ or $a_{1}\vee b_{2}=b_{1}$ and $a_{1}\wedge b_{2}=a_{2}$.  In that case we write $(a_{1},b_{1})\perspective(a_{2},b_{2})$.  This is illustrated in Figure~\ref{fig:perspective_covers}.

For $a,b\in L$ with $a\leq b$, we define the associated \defn{interval} by
\begin{displaymath}
    [a,b]\defs\{c\in L\colon a\leq c\leq b\}.
\end{displaymath}

A subset $C\subseteq L$ is a \defn{chain} if it can be written as $C=\{c_{0},c_{1},\ldots,c_{k}\}$ such that $c_{0}<c_{1}<\cdots<c_{k}$.  The \defn{length} of such a chain is $k$.  A chain is \defn{maximal} if $c_{0}=\least$, $c_{k}=\grtst$ and $(c_{i-1},c_{i})\in\Covers(\Lattice)$ for all $i\in[k]$.  The \defn{length} of $\Lattice$, denoted by $\len(\Lattice)$ is the maximum length of a maximal chain of $\Lattice$.  

An element $j\in L\setminus\{\least\}$ is \defn{join irreducible} if $j=a\vee b$ implies $j\in\{a,b\}$.  Since $\Lattice$ is finite, $j\in L$ is join irreducible if and only if $j$ covers a unique element, denoted by $j_{*}$.  We denote the set of all join-irreducible elements of $\Lattice$ by $\JI(\Lattice)$.  If $j_{*}=\least$, then $j$ is an \defn{atom}.  

\begin{lemma}\label{lem:covers_perspective_irreducible}
	For every $(a,b)\in\Covers(\Lattice)$, there exists $j\in\JI(\Lattice)$ such that $(a,b)\perspective(j_{*},j)$.
\end{lemma}
\begin{proof}
	We proceed by induction on $\len(\Lattice)$.  If $\len(\Lattice)=1$, then $L=\{\least,\grtst\}$, $\Covers(\Lattice)=\bigl\{(\least,\grtst)\bigr\}$ and $\JI(\Lattice)=\{\grtst\}$.  The claim then holds trivially, because $(\least,\grtst)\perspective(\least,\grtst)$.  
	
	Now suppose that $\len(\Lattice)>1$ and that the claim holds for all lattices of length strictly smaller than $\len(\Lattice)$.  Pick $(a,b)\in\Covers(\Lattice)$.  
	
	(i) If $b\neq\grtst$, then by induction we can find $j\in\JI\bigl([\least,b]\bigr)$ such that $(a,b)\perspective(j_{*},j)$, because $\len\bigl([\least,b]\bigr)<\len(\Lattice)$.  But then, $j\in\JI(\Lattice)$ because $[\least,b]$ is an interval of $\Lattice$ and it holds $(a,b)\perspective(j_{*},j)$ in $\Lattice$ as well.
	
	(ii) If $b=\grtst$, then there are two options.  Either $\grtst\in\JI(\Lattice)$ or not.  In the first case, we have $a=\grtst_{*}$ and the claim holds trivially, because $(a,b)\perspective(\grtst_{*},\grtst)$.  Otherwise, there exists $c\in L$ such that $c\neq a$ and $c\lessdot \grtst$.  Let $z=a\wedge c$.  Since $\Lattice$ is finite, there must exist $d\in L$ such that $z\lessdot d\leq c\lessdot b$ and $d\not\leq a$.  Since $\len\bigl([\least,d])<\len(\Lattice)$, by induction we can find $j\in\JI\bigl([\least,d]\bigr)$ such that $(z,d)\perspective(j_{*},j)$.  But then $j\in\JI(\Lattice)$, because $[\least,d]$ is an interval of $\Lattice$ and it holds $(z,d)\perspective(j_{*},j)$ in $\Lattice$ as well.  By definition, we have $j\leq d\leq c\lessdot b$.  If $j\leq a$, then we have $j\leq a\wedge c=z$ which contradicts $j\vee z=d$. Thus, $j\not\leq a$ and we conclude $j\vee a=b$.  Moreover, $j\wedge a\leq j_{*}$, because $j\not\leq a$ and $j\in\JI(\Lattice)$, and $j_{*}=j\wedge z\leq j\wedge a$, because $z\leq a$.  It follows that $j_{*}=j\wedge a$, and we conclude $(a,b)\perspective(j_{*},j)$.
\end{proof}

Inspired by this property, we consider for any maximal chain $C=\{c_{0},c_{1},\ldots,c_{k}\}$ the map
\begin{equation}
	\gammairr_{C}\colon\JI(\Lattice)\to[k], \quad j\mapsto\min\{s\colon j\leq c_{s}\}.
\end{equation}

Note that it is not necessarily the case that any join-irreducible element is \emph{perspective} to some covering pair in $C$; see for instance Figure~\ref{fig:perspective_covers}.  As a corollary, we obtain the following well-known relation between the length of $\Lattice$ and the number of join-irreducibles.

\begin{corollary}\label{cor:gamma_surjective}
	For any maximal chain $C$, the map $\gammairr_{C}$ is surjective.  Consequently, $\len(\Lattice)\leq\bigl\lvert\JI(\Lattice)\bigr\rvert$.
\end{corollary}
\begin{proof}
	Let $C=\{c_{0},c_{1},\ldots,c_{k}\}$ be a maximal chain of $\Lattice$.  For $s\in[k]$, Lemma~\ref{lem:covers_perspective_irreducible} implies that there exists $j_{s}\in\JI(\Lattice)$ with $(c_{s-1},c_{s})\perspective({j_{s}}_{*},j_{s})$.  By definition, we have $j_{s}\leq c_{s}$ and $j_{s}\not\leq c_{s-1}$ which implies $\gammairr_{C}(j_{s})=s$.  This yields the first part of the statement.  The second part of the statement follows if we take $k=\len(\Lattice)$.
\end{proof}

\subsection{Some generalizations of distributive lattices}

\subsubsection{Join-extremal lattices} 

Following \cite{markowsky92primes}, we award the lattices which satisfy equality in Corollary~\ref{cor:gamma_surjective} a special name.  The lattice $\Lattice$ is \defn{join extremal} if $\len(\Lattice)=\bigl\lvert\JI(\Lattice)\bigr\rvert$.  It is \defn{extremal} if both $\Lattice$ and $\Lattice^{\dual}$ are join extremal.

\begin{corollary}\label{cor:gamma_bijective}
	Let $C$ be a maximal chain of $\Lattice$ with $\len(\Lattice)=\lvert C\rvert-1$.  The map $\gamma_{C}$ is a bijection if and only if $\Lattice$ is join extremal.
\end{corollary}
\begin{proof}
	Let $C=\{c_{0},c_{1},\ldots,c_{k}\}$ with $k=\len(\Lattice)$.  If $\Lattice$ is join extremal, then $k=\bigl\lvert\JI(\Lattice)\bigr\rvert$, and $\gammairr_{C}$ is a surjective map between equinumerous sets.  Therefore it must be a bijection.  If $\Lattice$ is not join extremal, then $k<\bigl\lvert\JI(\Lattice)\bigr\rvert$, and the pigeon-hole principle tells us that $\gammairr_{C}$ cannot be injective.
\end{proof}

\subsubsection{Left-modular lattices}

Let $a,b,c\in L$ with $b<c$.  Then, these elements satisfy the \defn{modular inequality} 
\begin{equation}\label{eq:mod}
	(b\vee a)\wedge c \geq b\vee(a\wedge c).
\end{equation}
If this holds with equality for all $b<c$, then the element $a$ is \defn{left modular}.  Following \cite{blass97mobius}, the lattice $\Lattice$ is \defn{left modular} if it has a maximal chain of size $\len(\Lattice)+1$ which consists entirely of left-modular elements.
%
See also \cites{liu00left,thomas06analogue,thomas19rowmotion} for more background.

\subsubsection{Join-semidistributive lattices}

Lastly, $\Lattice$ is \defn{join semidistributive} if for all $a,b,c\in L$ it holds that
\begin{equation}\label{eq:jsd}
	a\vee b=a\vee c\quad\text{implies}\quad a\vee b=a\vee(b\wedge c).
\end{equation}
If $\Lattice$ and $\Lattice^{\dual}$ are join semidistributive, then $\Lattice$ is \defn{semidistributive}.

\section{Proof of the Main Result}

\begin{figure}
	\centering
	\begin{subfigure}[t]{.27\textwidth}
		\centering
		\begin{tikzpicture}
			\def\x{1};
			\draw(2*\x,1*\x) node[draw,circle,scale=.67](a1){};
			\draw(1*\x,2*\x) node[draw,circle,scale=.67](a2){};
			\draw(3*\x,2*\x) node[draw,circle,scale=.67](a3){};
			\draw(1*\x,3*\x) node[draw,circle,scale=.67](a4){};
			\draw(3*\x,3*\x) node[draw,circle,scale=.67](a5){};
			\draw(2*\x,4*\x) node[draw,circle,scale=.67](a6){};
			\draw[thick](a1) -- (a2);
			\draw[thick](a1) -- (a3);
			\draw[thick](a2) -- (a4);
			\draw[thick](a3) -- (a5);
			\draw[thick](a4) -- (a6);
			\draw[thick](a5) -- (a6);
		\end{tikzpicture}
		\caption{A join-semidistributive lattice that is neither join extremal nor left modular.}
		\label{fig:hexagon}
	\end{subfigure}
	\hspace*{.5cm}
	\begin{subfigure}[t]{.27\textwidth}
		\centering
		\begin{tikzpicture}
			\def\x{1};
			\draw(2*\x,1*\x) node[draw,circle,scale=.67](a1){};
			\draw(1*\x,2*\x) node[draw,circle,scale=.67](a2){};
			\draw(2*\x,2*\x) node[draw,circle,scale=.67](a3){};
			\draw(3*\x,2*\x) node[draw,circle,scale=.67](a4){};
			\draw(2*\x,3*\x) node[draw,circle,scale=.67](a5){};
			\draw[thick](a1) -- (a2);
			\draw[thick](a1) -- (a3);
			\draw[thick](a1) -- (a4);
			\draw[thick](a2) -- (a5);
			\draw[thick](a3) -- (a5);
			\draw[thick](a4) -- (a5);
		\end{tikzpicture}
		\caption{A left-modular lattice that is neither join extremal nor join semidistributive.}
		\label{fig:m3}
	\end{subfigure}
	\hspace*{.5cm}
	\begin{subfigure}[t]{.27\textwidth}
		\centering
		\begin{tikzpicture}
			\def\x{1};
			\draw(2*\x,1*\x) node[draw,circle,scale=.67](a1){};
			\draw(1*\x,2*\x) node[draw,circle,scale=.67](a2){};
			\draw(2*\x,2*\x) node[draw,circle,scale=.67](a3){};
			\draw(3*\x,2*\x) node[draw,circle,scale=.67](a4){};
			\draw(1*\x,3*\x) node[draw,circle,scale=.67](a5){};
			\draw(3*\x,3*\x) node[draw,circle,scale=.67](a6){};
			\draw(2.25*\x,2.75*\x) node[draw,circle,scale=.67](a7){};
			\draw(2*\x,3*\x) node[draw,circle,scale=.67](a8){};
			\draw(2*\x,4*\x) node[draw,circle,scale=.67](a9){};
			\draw[thick](a1) -- (a2);
			\draw[thick](a1) -- (a3);
			\draw[thick](a1) -- (a4);
			\draw[thick](a2) -- (a5);
			\draw[thick](a2) -- (a8);
			\draw[thick](a3) -- (a5);
			\draw[thick](a3) -- (a6);
			\draw[thick](a4) -- (a6);
			\draw[thick](a4) -- (a7);
			\draw[thick](a5) -- (a9);
			\draw[thick](a6) -- (a9);
			\draw[thick](a7) -- (a8);
			\draw[thick](a8) -- (a9);
		\end{tikzpicture}
		\caption{A join-extremal lattice that is neither left modular nor join semidistributive.}
		\label{fig:extremal}
	\end{subfigure}
	\caption{Some lattices.}
	\label{fig:lattices}
\end{figure}
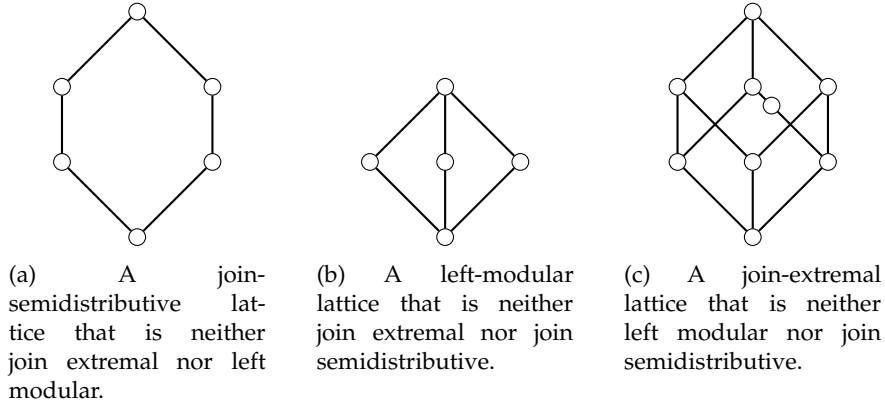

In general, there are no implications among the three lattice properties defined before.  Figure~\ref{fig:hexagon} shows a join-semidistributive lattice that is neither join extremal nor left modular; Figure~\ref{fig:m3} shows a left-modular lattice that is neither join extremal nor join semidistributive; Figure~\ref{fig:extremal} shows a join-extremal lattice that is neither join semidistributive nor left modular.

The main motivation for this article comes from the following result that relates extremality and left-modularity for semidistributive lattices.

\begin{theorem}[\cite{thomas19rowmotion}*{Theorem~1.4}]\label{thm:extremal_is_trim}
	Every semidistributive, extremal lattice is left modular.
\end{theorem}

Our main contribution is the converse to Theorem~\ref{thm:extremal_is_trim} on the level of join-semidistributive lattices.  

\begin{theorem}\label{thm:semidistributive_left_modular_is_trim}
	Every join-semidistributive, left-modular lattice is join extremal.
\end{theorem}
\begin{proof}
	Suppose that $\Lattice$ is join semidistributive and left modular with $\len(\Lattice)=k$.  Fix a chain $C=\{c_{0},c_{1},\ldots,c_{k}\}$ of left-modular elements.  We proceed by way of contradiction and assume that $\gammairr_{C}$ is not injective.  Let $j_{1},j_{2}\in\JI(\Lattice)$ with $j_{1}\neq j_{2}$ and $\gammairr_{C}(j_{1})=s=\gammairr_{C}(j_{2})$.
	
	This implies, by construction, that $j_{i}\not\leq c_{s-1}$ for $i\in\{1,2\}$.  Since $c_{s-1}\lessdot c_{s}$ and $j_{i}\leq c_{s}$ it follows that $c_{s-1}<j_{i}$ if and only if $j_{i}=c_{s}$.  Thus, without loss of generality, we may assume that $j_{1}$ and $c_{s-1}$ are incomparable.
	We conclude:
	\begin{equation}\label{eq:A}
		j_{1}\vee c_{s-1} = c_{s} = j_{2}\vee c_{s-1}.
	\end{equation}
    Together with \eqref{eq:jsd}, we get
	\begin{equation}\label{eq:B}
		c_{s} = c_{s-1}\vee(j_{1}\wedge j_{2}).
	\end{equation}

	Now we show that $j_{1}$ and $j_{2}$ are incomparable, and assume that $j_{1}<j_{2}$.  Since $c_{s-1}$ is left modular, we obtain using \eqref{eq:A}:
	\begin{displaymath}
		j_{2} = c_{s}\wedge j_{2} = (j_{1}\vee c_{s-1})\wedge j_{2} = j_{1}\vee(c_{s-1}\wedge j_{2})\leq j_{1}\vee {j_{2}}_{*} \leq {j_{2}}_{*},
	\end{displaymath}
	which is a contradiction.  Switching the roles of $j_{1}$ and $j_{2}$ discards the case $j_{2}<j_{1}$.  Since $j_{1}\neq j_{2}$ by assumption, these elements must be incomparable.

	Using the left-modularity of $c_{s-1}$ once again, we obtain from \eqref{eq:B} that
	\begin{equation}\label{eq:C}
		j_{1} = c_{s}\wedge j_{1} = \bigl((j_{1}\wedge j_{2})\vee c_{s-1}\bigr)\wedge j_{1} = (j_{1}\wedge j_{2})\vee(c_{s-1}\wedge j_{1}).
	\end{equation}
	Since $j_{1}\in\JI(\Lattice)$, it follows that $j_{1}=j_{1}\wedge j_{2}$ or $j_{1}=c_{s-1}\wedge j_{1}$.  However, we have already established that $j_{1}$ and $j_{2}$ are incomparable (which discards the first case) and that $j_{1}$ and $c_{s-1}$ are incomparable (which discards the second case).  This contradiction shows that our assumption must have been wrong, meaning that $\gammairr_{C}$ is injective.  
	
	Together with Corollary~\ref{cor:gamma_surjective} we get that $\gammairr_{C}$ is bijective, so that Corollary~\ref{cor:gamma_bijective} implies that $\Lattice$ is join extremal.
\end{proof}

Since Theorem~\ref{thm:semidistributive_left_modular_is_trim} is stated for join-semidistributive lattices, it is natural to wonder if we can weaken any of the assumptions of Theorem~\ref{thm:extremal_is_trim}.  For keeping the notation simple, we define $\MI(\Lattice)\defs\JI(\Lattice^{\dual})$.  First of all, we record the following lemma.

\begin{lemma}[\cite{freese95free}*{Corollary~2.55}]\label{lem:semidistributive_extremal}
	If $\Lattice$ is semidistributive, then $\bigl\lvert\JI(\Lattice)\bigr\rvert=\bigl\lvert\MI(\Lattice)\bigr\rvert$.  
\end{lemma}

In particular, semidistributive, join-extremal lattices are extremal.  This gives us the following corollary.

\begin{corollary}\label{cor:left_modular_is_extremal}
	For semidistributive lattices, left-modularity and extremality are equivalent.
\end{corollary}

For $j\in\JI(\Lattice)$, let us consider the set
\begin{displaymath}
	K(j) \defs \bigl\{a\in L\colon j_{*}\leq a\;\text{but}\;j\not\leq a\bigr\}.
\end{displaymath}
If $K(j)$ has a unique maximal element, then we write $\kappa(j)$ for this element.  In other words $\kappa(j)\defs\bigvee K(j)$ if this join exists.  

\begin{lemma}\label{lem:kappa_maximal}
	Every maximal element of $K(j)$ is in $\MI(\Lattice)$.  Consequently, $\kappa(j)\in\MI(\Lattice)$ if it exists.
\end{lemma}
\begin{proof}
	Let $j\in\JI(\Lattice)$.  Since $j_{*}\in K(j)$, the set $K(j)$ is not empty.  Let $m\in K(j)$ be maximal.  Assume that $m\notin\MI(\Lattice)$.  If $m=\grtst$, then $j\leq m$, which is a contradiction.  Therefore, $m\neq\grtst$.  By construction, this means that there exist (at least) two distinct elements $a_{1},a_{2}\in L$ such that $m\lessdot a_{1}$ and $m\lessdot a_{2}$.  Since $m\in K(j)$ we conclude that $j_{*}\leq m$ but $j\not\leq m$.  In particular $j_{*}\leq a_{1}$ and $j_{*}\leq a_{2}$.  Since $m$ is maximal, we conclude that $a_{1},a_{2}\notin K(j)$ meaning that necessarily $j\leq a_{1}$ and $j\leq a_{2}$.  But then $j\leq a_{1}\wedge a_{2}=m$, a contradiction.
\end{proof}

\begin{lemma}\label{lem:kappa_disjoint}
	Let $\Lattice$ be join semidistributive.  If $j,j'\in\JI(\Lattice)$ are distinct, then there does not exist $m\in L$ which is maximal in both $K(j)$ and $K(j')$.
\end{lemma}
\begin{proof}
	Assume that there exists $m\in L$ which is a maximal element of both $K(j)$ and $K(j')$.  By Lemma~\ref{lem:kappa_maximal}, $m\in\MI(\Lattice)$.  Let $m^{*}$ denote the unique element in $\Lattice$ that covers $m$.
	
	By assumption, $m^{*}\notin K(j)\cup K(j')$, which means that $j\leq m^{*}$ and $j'\leq m^{*}$.  Then, $j\vee m=m^{*}$ and $j'\vee m=m^{*}$, because $m\lessdot m^{*}$.  If $z=j\wedge j'$, then \eqref{eq:jsd} implies that $m^{*}=m\vee z$.  Since $j\in\JI(\Lattice)$ and $j\neq j'$, it follows that $z\leq j_{*}\leq m$, which implies $m\vee z=m$.  We thus obtain the contradiction $m^{*}=m$.
\end{proof}

\begin{lemma}\label{lem:extremal_semidistributive}
	Every join-semidistributive lattice with $\bigl\lvert\JI(\Lattice)\bigr\rvert=\bigl\lvert\MI(\Lattice)\bigr\rvert$ is semidistributive.
\end{lemma}
\begin{proof}
	For $j\in\JI(\Lattice)$ we denote by $M(j)$ the set of maximal elements of $K(j)$.  By Lemma~\ref{lem:kappa_maximal}, $M(j)\subseteq\MI(\Lattice)$.  Moreover, by Lemma~\ref{lem:kappa_disjoint}, if $j,j'\in\JI(\Lattice)$ with $j\neq j'$ we have $M(j)\cap M(j')=\emptyset$.  Let $\JI(\Lattice)=\{j_{1},j_{2},\ldots,j_{k}\}$ and write $m_{\ell}\defs\bigl\lvert M(j_{\ell})\bigr\rvert$ for $\ell\in[k]$.  Since $K(j)\neq\emptyset$ by construction it follows that $m_{\ell}\geq 1$ for all $\ell\in[k]$.
	
	If $\bigl\lvert\MI(\Lattice)\bigr\rvert=k$, then we obtain
	\begin{displaymath}
		k \leq m_{1}+m_{2}+\cdots+m_{k} \leq \bigl\lvert\MI(\Lattice)\bigr\rvert = k,
	\end{displaymath}
	which enforces $m_{\ell}=1$ for all $\ell\in[k]$.  Consequently, $\kappa(j)$ exists for all $j\in\JI(\Lattice)$.  By \cite{freese95free}*{Theorem~2.56} it follows that $\Lattice^{\dual}$ is join-semidistributive.  
\end{proof}

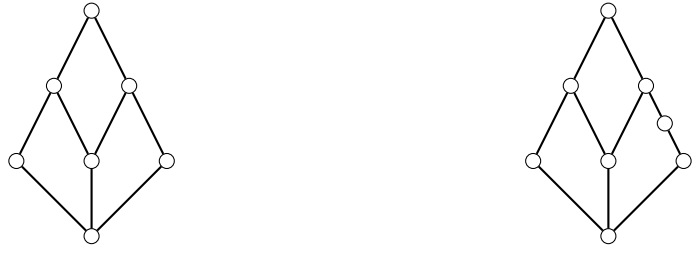
\begin{figure}
	\centering
		\begin{subfigure}[t]{.45\textwidth}
		\centering
		\begin{tikzpicture}\small
			\def\x{1};
			\draw(2*\x,1*\x) node[draw,circle,scale=.67](a1){};
			\draw(1*\x,2*\x) node[draw,circle,scale=.67](a2){};
			\draw(2*\x,2*\x) node[draw,circle,scale=.67](a3){};
			\draw(3*\x,2*\x) node[draw,circle,scale=.67](a4){};
			\draw(1.5*\x,3*\x) node[draw,circle,scale=.67](a5){};
			\draw(2.5*\x,3*\x) node[draw,circle,scale=.67](a6){};
			\draw(2*\x,4*\x) node[draw,circle,scale=.67](a7){};
			\draw[thick](a1) -- (a2);
			\draw[thick](a1) -- (a3);
			\draw[thick](a1) -- (a4);
			\draw[thick](a2) -- (a5);
			\draw[thick](a3) -- (a5);
			\draw[thick](a3) -- (a6);
			\draw[thick](a4) -- (a6);
			\draw[thick](a5) -- (a7);
			\draw[thick](a6) -- (a7);
		\end{tikzpicture}
		\caption{A join-semidistributive, join-extremal, left-modular lattice that is not semidistributive.}
		\label{fig:left_modular_not_semidistributive}
	\end{subfigure}
	\hspace*{1cm}
	\begin{subfigure}[t]{.45\textwidth}
		\centering
		\begin{tikzpicture}\small
			\def\x{1};
			\draw(2*\x,1*\x) node[draw,circle,scale=.67](a1){};
			\draw(1*\x,2*\x) node[draw,circle,scale=.67](a2){};
			\draw(2*\x,2*\x) node[draw,circle,scale=.67](a3){};
			\draw(3*\x,2*\x) node[draw,circle,scale=.67](a4){};
			\draw(2.75*\x,2.5*\x) node[draw,circle,scale=.67](a5){};
			\draw(1.5*\x,3*\x) node[draw,circle,scale=.67](a6){};
			\draw(2.5*\x,3*\x) node[draw,circle,scale=.67](a7){};
			\draw(2*\x,4*\x) node[draw,circle,scale=.67](a8){};
			\draw[thick](a1) -- (a2);
			\draw[thick](a1) -- (a3);
			\draw[thick](a1) -- (a4);
			\draw[thick](a2) -- (a6);
			\draw[thick](a3) -- (a6);
			\draw[thick](a3) -- (a7);
			\draw[thick](a4) -- (a5);
			\draw[thick](a5) -- (a7);
			\draw[thick](a6) -- (a8);
			\draw[thick](a7) -- (a8);
		\end{tikzpicture}
		\caption{A join-semidistributive, join-extremal lattice that is not left modular.}
		\label{fig:join_extremal_not_left_modular}
	\end{subfigure}
	\caption{Some join-semidistributive lattices.}
	\label{fig:jsd_lattices}
\end{figure}

As a consequence, every join-semidistributive, extremal lattice is semidistributive.  Finally, Figure~\ref{fig:join_extremal_not_left_modular} shows a join-semidistributive, join-extremal lattice that is not left modular.  Therefore, we cannot weaken the assumptions of Theorem~\ref{thm:extremal_is_trim} in order to guarantee left-modularity.

\section{Shellability}
	\label{sec:shellability}
We now briefly touch a topological aspect of the lattices that we consider.  An \defn{edge labeling} of $\Lattice$ is any map
\begin{displaymath}
	\lambda\colon\Covers(\Lattice)\to M
\end{displaymath}
for some set $M$.  If $C=\{c_{0},c_{1},\ldots,c_{k}\}$ is a maximal chain of $\Lattice$, then 
\begin{displaymath}
	\lambda(C) \defs \bigl(\lambda(c_{0},c_{1}),\lambda(c_{1},c_{2}),\ldots,\lambda(c_{k-1},c_{k})\bigr)
\end{displaymath}
is the associated label vector.  If $M\subseteq\mathbb{N}$, then $C$ is \defn{increasing} if $\lambda(C)$ is strictly increasing.  Following \cites{bjorner80shellable,bjorner96shellable}, the edge labeling $\lambda$ is an \defn{EL-labeling} if every interval $[a,b]$ of $\Lattice$ contains a unique increasing maximal chain, and this maximal chain has the lexicographically smallest label vector among all maximal chains in $[a,b]$.  The lattice $\Lattice$ is \defn{EL-shellable} if it admits an EL-labeling.  EL-shellable lattices have remarkable properties, for instance the order complex of $\Lattice\setminus\{\least,\grtst\}$ is a shellable, hence Cohen--Macaulay, complex.

If $\Lattice$ is left modular with left-modular chain $C=\{c_{0},c_{1},\ldots,c_{k}\}$, then the map $\gammairr_{C}$ induces an edge labeling of $\Lattice$ by setting
\begin{displaymath}
	\lmlabeling{C}\colon\Covers(\Lattice)\to[k], \quad (a,b)\mapsto\min\bigl\{\gammairr_{C}(j)\colon j\in\JI(\Lattice), a\vee j=b\bigr\}.
\end{displaymath}
It was shown in \cite{liu99left} that $\lmlabeling{C}$ is an EL-labeling, which yields the following result.

\begin{proposition}[\cite{liu99left}]\label{prop:left_modular_shellable}
	Every left-modular lattice is EL-shellable.
\end{proposition}

The converse of Proposition~\ref{prop:left_modular_shellable} is not true, see for instance the lattice in Figure~\ref{fig:extremal}.  Together with Theorem~\ref{thm:extremal_is_trim}, we obtain the following.

\begin{corollary}
	Every semidistributive, extremal lattice is EL-shellable.
\end{corollary}

Theorems~\ref{thm:extremal_is_trim} and \ref{thm:semidistributive_left_modular_is_trim} state that we cannot distinguish between extremality and left-modularity on the level of semidistributive lattices.  What about shellability?  Proposition~\ref{prop:left_modular_shellable} on the level of join-semidistributive lattices reads as follows.

\begin{corollary}\label{cor:left_modular_shellable}
	Every join-semidistributive, left-modular lattice is EL-shellable.
\end{corollary}

The converse does not hold by virtue of the lattice in Figure~\ref{fig:join_extremal_not_left_modular}.  This lattice is join semidistributive and EL-shellable, but not left modular.  But are we able to distinguish left-modularity and EL-shellability on the level of semidistributive lattices?

\begin{question}\label{qu:shellable_left_modular}
	Is every semidistributive, EL-shellable lattice necessarily left modular?
\end{question}

We note that in view of Theorem~\ref{thm:semidistributive_left_modular_is_trim}, if Question~\ref{qu:shellable_left_modular} is true, then any semidistributive, EL-shellable lattice must necessarily be extremal.  Is that always the case?

\begin{question}\label{qu:shellable_extremal}
    Is every semidistributive, EL-shellable lattice necessarily extremal?
\end{question}

Figure~\ref{fig:m3} shows an EL-shellable lattice that is not extremal (and it is also not semidistributive).  

We summarize the implications among the various families of lattices in Figure~\ref{fig:big_picture}.  In this figure, the gray boxes denote the classes of lattices having the properties stated in the corresponding row and column headers.  For instance, the green arrow from second box in the bottom row to the third box in the bottom row represents the statement:
\begin{center}
    ``Every semidistributive, join-extremal lattice is left modular.''
\end{center}
which is true by virtue of Theorem~\ref{thm:extremal_is_trim} and Lemma~\ref{lem:semidistributive_extremal}.


\end{document}